\def\softd{{\leavevmode\setbox1=\hbox{d}%
\hbox to 1.05\wd1{d\kern-0.4ex{\char039}\hss}}}
\def\softt{{\leavevmode\setbox1=\hbox{t}%
\hbox to \wd1{t\kern-0.6ex{\char039}\hss}}}
\def\softl{l\kern-0.45ex\raise0.1ex\hbox{'}\kern-0.10ex}
\def\softL{L\kern-0.8ex\raise0.1ex\hbox{'}\kern0.1ex}
\documentstyle[11pt,A4]{report}
\newcounter{newexample}
\newcommand{\example}{\vspace*{1.5ex} \noindent
                      \protect\refstepcounter{newexample}
                      \underline{{\bf Example \thenewexample.}}}

\title
{\vspace*{-5cm} {\sc Slovak Academy of Sciences \\
 Institute of Experimental Physics} \\
\vspace*{4cm}
\sc The Laplace Transform Method \\
            for Linear Differential Equations \\
            of the Fractional Order}
\author{{\bf Igor Podlubny} \protect \\
Department of Control Engineering \\
Faculty of Mining, University of Technology\\
B.Nemcovej 3, 042 00 Kosice, Slovakia \\
Phone: (+42 95) 39772 \\
Fax: (+42 95) 36618 \\
E-mail: podlbn@ccsun.tuke.sk}
\date{\vspace*{2cm} UEF-02-94, \hspace*{3em} June 1994 \\[1cm]
        \parbox[t]{12cm}{\footnotesize
        \rule{12cm}{0.5pt} \\
        \indent
        This is the e-Print version of the pre-print printed in 1994.\\
        Misprints were corrected.\\[1ex]
        \indent
        Present address: \\
        Department of Management and Control
        Engineering, Faculty of B.E.R.G., Technical University of Kosice,
        B.Nemcovej 3, 042 00 Kosice, Slovak Republic.
        Phone: (+421 95) 6339772; fax: (+421 95) 6336618;
        e-mail: as above.}
}

\begin{document}
\begin{sloppypar}
\maketitle

\newpage\thispagestyle{empty}
\copyright\ 1994, RNDr.~Igor~Podlubny,~CSc.

\vspace*{20cm}
\noindent
This publication was typeset by \LaTeX.

\setcounter{page}{0}
\tableofcontents

\newpage\thispagestyle{empty}
\mbox{~}

\chapter*{Preface}
\addcontentsline{toc}{chapter}{Preface}

Differential equations of the fractional order appear more and more
frequently in different research areas and engineering applications.
An effective and easy-to-use method for solving
such equations is needed.

However, known methods have certain disadvantages.
Methods, described in details in \cite{OS,Miller-Ross,Bagley-Calico}
for fractional differential equations
of the rational order, do not work in the case of an arbitrary
real order. On the other hand, there is an iteration method
described in \cite{SKM}, which allows solution of fractional
differential equations of an arbitrary real order, but
it works effectively only for relatively simple equations,
as well as the series method \cite{OS,Friedrich-1991}.
Other authors (e.g. \cite{Bagley-Calico,Caputo-Mainardi}) used
in their investigations the one-parameter Mittag-Leffler function
$E_{\alpha}(z)=\sum_{k=0}^{\infty} \frac{z^k}{\Gamma (\alpha k + 1)}$.
Still other authors \cite{Schneider-Wyss,GN}
prefer the Fox $H$-function \cite{Fox}, which seems to be too general
to be frequently used in applications.

Instead of this variety of different methods, we introduce a method
which is free of the mentioned disadvantages and suitable for
a wide class of initial value problems for fractional differential equations.
The method uses the Laplace transform technique and is based on the
formula of the Laplace transform of the Mittag-Leffler function
in two parameters $E_{\alpha, \beta}(z)$.
We hope that the described method could be useful for obtaining solutions
of different applied problems, appearing in physics, chemistry,
electrochemistry, engineering, financing and banking, etc.
To outline the area of the method's applicability, we have
included in the bibliography also the works by different authors
\cite{Caputo-1967}--\cite{Axtell-1990},
in which fractional linear differential equations appears or
which could serve as a basis for obtaining such equations.

\vspace{1.5ex}

This work deals with solution of the fractional linear
differential equations with constant coefficients
and consists of four short chapters.

In Chapter 1 we present some auxillary tools which are necessary
for using the method. The reader can find the definition there
and some important properties of the Mittag-Leffler function
in two parameters and the Wright function.
        The basic result, presented in Chapter 1, is the Laplace
transform of the Mittag-Leffler function and its derivatives.
        Besides that, we introduce two tools necessary for testing
candidate solutions by direct substitution in corresponding equations:
fractional derivatives of the Mittag-Leffler function and
a rule for the fractional differentiation of integrals
depending on a parameter.

In Chapter 2 we give solutions to some initial-value problems
for ''standard'' fractional differential equations.
Some of them were solved by other authors earlier by other methods,
and the comparison in such cases just underline the simplicity
and the power of our approach.

In Chapter 3 we extend the proposed method for the case of so-called
''sequential'' fractional differential equations (we adopted the convenient
terminology of Miller and Ross). For this purpose, we obtained
the Laplace transform for the ''sequential'' fractional derivative.
The ''sequential'' analogues of the problems, solved in Chapter 2,
are considered. Naturally, we arrive at solutions which are different
from those obtained in Chapter 2.

However, there is something common in solutions of the corresponding
''standard'' and ''sequential'' fractional differential equations:
they both have the same {\em fractional Green's function}. In Chapter 4 we
give our definition of the fractional Green's function and some
of its properties, necessary for constructing solutions of
initial-value problems for fractional linear differential equations
with constant coefficients.

We give the solution of the initial-value problem for
the ordinary fractional linear differential
equation with constant coefficients using only its Green's function.
Due to this result, the solution of such initial value problems
reduces to finding the fractional Green's functions.
We obtained the explicit expressions for the fractional
Green's function for some special cases (one- , two- , three- and
four-term equations).

The explicit expression for an arbitrary
fractional linear ordinary differential equation with constant
coefficients ends this work.

\vspace{1.5ex}

In what follows below, $_{a}D_{t}^{\alpha}$ means the Riemann-Liouville
fractional derivative of order $\alpha$
(e.g., \cite{OS,Miller-Ross}):
$$
   _{a}D_{t}^{\alpha}f(t)=
    \frac{1}{\Gamma (n -\alpha)}
    \frac{d^{n}}{dt^{n}}
    \int_{a}^{t}
    \frac{f(\tau)}{(t-\tau)^{\alpha - n + 1}}d\tau ,
    \hspace{2em}
    (n-1 < \alpha <n)
$$

\vspace{1ex}

The author thanks Ms.~Serena~Yeo for checking the language of this work.

\cleardoublepage
\chapter{Introduction}

\section{Mittag-Leffler function in two parameters}

The function, which plays a very important role in this work,
was in fact introduced by Agarwal \cite{Agarwal}. A number
of relationships for this function were obtained by Agarwal
and Humbert \cite{Humbert-Agarwal} using the Laplace transform technique.
This function could have been called the Agarwal function.
However, Agarwal and Humbert generously
left the same notation as for the one-parameter
Mittag-Leffler function, and that was the reason that
now the two-parameter function is called the Mittag-Leffler function.
We will use the name and the notation used
in the fundamental handbook on special functions \cite{HTF-3}.
In spite of using the same notation as Agarwal,
the definition given there differs from Agarwal's definition
by a non-constant factor.

\vspace{1.5ex}
\noindent
\underline{{\bf Definition.}}
A two-parameter function of the Mittag-Leffler type is defined
by the series expansion \cite{HTF-3}
\begin{equation} \label{ML-Definition}
E_{\alpha, \beta}(z) = \sum_{k=0}^{\infty}\frac{z^k}{\Gamma (\alpha k + \beta)},
\hspace{3em} (\alpha > 0, \hspace{1em} \beta > 0)
\end{equation}

\noindent
\underline{{\bf Relation to other functions.}}
There are some relationships given in \cite[\S 18.1]{HTF-3}:
\begin{eqnarray}
 & & E_{1,1}(z) = e^{z}, \hspace{2em} E_{1,2}(z) = \frac{\mbox{$e^z -1$}}{z},
                  \label{Exp}\\
 & & E_{2,1}(z) = \cosh (\sqrt{z}), \hspace{2em}
                E_{2,2}(z) = \frac{\mbox{$\sinh (\sqrt{z})$}}{\sqrt{z}}, \\
 & & E_{1/2, 1}(\sqrt{z}) = \frac{2}{\sqrt{\pi}} e^{-z} \mbox{erfc}(-\sqrt{z}).
\end{eqnarray}

For $\beta = 1$ we obtain the Mittag-Leffler function in one parameter:
\begin{equation}
E_{\alpha, 1}(z) = \sum_{k=0}^{\infty}\frac{z^k}{\Gamma (\alpha k + 1)}
                   \equiv E_{\alpha}(z).
\end{equation}
The function ${\cal E}_{t} (\nu,a)$, introduced in \cite{Miller-Ross}
for solving differential equations of the rational order,
is the particular case of the Mittag-Leffler function (\ref{ML-Definition}):
\begin{equation} \label{Miller-Ross-Function}
{\cal E}_{t}(\nu, a) = t^{\nu} \sum_{k=0}^{\infty}\frac{(at)^{k}}
                                                       {\Gamma (\nu + k + 1)}
                     = t^{\nu} E_{1, \nu + 1} (at).
\end{equation}

Rabotnov's \cite{Rabotnov} function $\ni \!\! _{\alpha}(\beta, t)$
is the particular case
of the Mittag-Leffler function (\ref{ML-Definition}) too:
\begin{equation} \label{Rabotnov-Function}
      \ni\ \!\!\! _{\alpha}(\beta, t) =
          t^{\alpha} \sum_{k=0}^{\infty}
                     \frac{\beta^{k}t^{k(\alpha +1)}}{\Gamma ((k+1)(1+\alpha))}
        = t^{\alpha} E_{\alpha+1, \alpha+1}(\beta t^{\alpha +1}).
\end{equation}

It follows from the relationships (\ref{Miller-Ross-Function})
and (\ref{Rabotnov-Function}) that the properties of the Miller-Ross
function and Rabotnov's function follows from the properties
of the Mittag-Leffler function in two parameters (\ref{ML-Definition}).

Plotnikov \cite[cf. \cite{Tseytlin}]{Plotnikov}
and Tseytlin \cite{Tseytlin} used in their investigations
two functions $Sc_{\alpha}(z)$ and $Cs_\alpha (z)$,
which they call the fractional sine and cosine.
Those functions are also just the particular cases of the
Mittag-Leffler function in two parameters:
\begin{equation}
      Sc_\alpha (z) = \sum_{n=0}^{\infty}
                        \frac{(-1)^n z^{(2-\alpha)n+1}}
                             {\Gamma ((2-\alpha)n+2)}
                    = z E_{2-\alpha, 2}(-z^{2-\alpha})
\end{equation}
\begin{equation}
      Cs_\alpha (z) = \sum_{n=0}^{\infty}
                        \frac{(-1)^n z^{(2-\alpha)n}}
                             {\Gamma ((2-\alpha)n+1)}
                    = z E_{2-\alpha, 2}(-z^{2-\alpha})
\end{equation}

Of course, properties of the fractional sine and cosine
follow from the properties of the Mittag-leffler
function (\ref{ML-Definition}).

\section{The Laplace transform of the \protect \\
          Mittag-Leffler function in two parameters}

As follows from relationship (\ref{Exp}),
the Mittag-Leffler function $E_{\alpha, \beta}(z)$ is a generalization
of the exponential function $e^{z}$ and the exponential function is
a particular case of the Mittag-Leffler function.

We will outline here the way to obtain the Laplace transform
of the Mittag-Leffler function with the help of the analogy
between this function and the function $e^{z}$.
For this purpose, let us obtain the Laplace transform of
the function $t^{k}e^{at}$ in an untraditional way.

First, let us prove that
\begin{equation} \label{Exp-Integral}
      \int_{0}^{\infty} e^{-t} e^{\pm zt} dt = \frac{1}{1 \mp z},
      \hspace{1em}
      \left| z \right| <1.
\end{equation}

Indeed, using the series expansion for $e^{z}$, we obtain
\begin{equation}
        \int_{0}^{\infty} e^{-t} e^{zt} dt = \frac{1}{1-z}
        = \sum_{k=0}^{\infty} \frac{(\pm z)^{k}}{k!}
                            \int_{0}^{\infty} e^{-t}t^{k} dt
        = \sum_{k=0}^{\infty} (\pm z)^{k} = \frac{1}{1 \mp z}.
\end{equation}

Second, we differentiate both sides of equation (\ref{Exp-Integral})
with respect to $z$. The result is
\begin{equation}
      \int_{0}^{\infty} e^{-t} t^{k} e^{\pm zt} dt = \frac{k!}{(1-z)^{k+1}},
      \hspace{1em}
      (|z| < 1),
\end{equation}
and after obvious substitutions we obtain
the well-known pair of the Laplace transform
of the function $t^{k}e^{\pm at}$:
\begin{equation}
    \int_{0}^{\infty} e^{-pt} t^{k} e^{\pm at} dt = \frac{k!}{(p \mp a)^{k+1}},
    \hspace{1em}
    (Re(p) > |a|).
\end{equation}

Let us now consider the Mittag-Leffler function (\ref{ML-Definition}).
Substitution of (\ref{ML-Definition}) in the integral below leads to
\begin{equation}\label{Basic}
\int_{0}^{\infty} e^{-t}t^{\beta-1}E_{\alpha, \beta}(zt^{\alpha}) dt
  = \frac{1}{1-z}, \hspace{2em} (\left| z \right| < 1),
\end{equation}
and we obtain from (\ref{Basic}) a pair
of the Laplace transforms of the function
$t^{\alpha k + \beta-1}E_{\alpha, \beta}^{(k)}(\pm zt^{\alpha})$, \hspace{0.3em}
($E_{\alpha,\beta}^{(k)}(y) \equiv \frac{d^{k}}{dy^{k}}E_{\alpha ,\beta}(y)$):
\begin{equation}\label{Laplace-Transform}
\int_{0}^{\infty} e^{-pt}t^{\alpha k +\beta-1}
              E_{\alpha, \beta}^{(k)}(\pm at^{\alpha}) dt
  = \frac{k! \,p^{\alpha - \beta}}{(p^{\alpha} \mp a)^{k+1}},
  \hspace{3em} ( Re(p) > \left| a \right|^{1/ \alpha}).
\end{equation}

\noindent
The particular case of (\ref{Laplace-Transform})
for $\alpha = \beta = \frac{1}{2}$
\begin{equation}\label{Laplace-Transform-Semi}
\int_{0}^{\infty} e^{-pt}t^{\frac{k-1}{2}}
        E_{\frac{1}{2}, \frac{1}{2}}^{(k)}(\pm a\sqrt{t}) dt
  = \frac{k!}{\left( \sqrt{p} \mp a \right)^{k+1}},
  \hspace{3em} ( Re(p) > a^{2} ).
\end{equation}
is useful for solving semidifferential equations considered
in \cite{OS,Miller-Ross}.

\section{The Wright function}

This function, related to the Mittag-Leffler function
in two parameters $E_{\alpha, \beta}(z)$,
was introduced by Wright
\cite[cf. \cite{HTF-3,Humbert-Agarwal}]{Wright-33}.
A number of useful relationships was obtained by Humbert and
Agarwal \cite{Humbert-Agarwal} with the help of
the Laplace transform.

For convenience we adopt here Mainardi's notation for the Wright function
$W(z; \alpha, \beta)$.

\noindent
\underline{{\bf Definition.}}
The Wright function is defined as \cite[formula 18.1(27)]{HTF-3}
\begin{equation} \label{W-Definition}
      W(z; \alpha, \beta) = \sum_{k=0}^{\infty} \frac{z^k}
                                   {k! \,\, \Gamma (\alpha k + \beta)}
\end{equation}

\noindent
\underline{{\bf Integral representation.}}
This function can be represented by the following integral
\cite[formula 18.1(29)]{HTF-3}:
\begin{equation} \label{W-integral-representation}
      W(z; \alpha, \beta) = \frac{1}{2 \pi i}
                            \int_{Ha} u^{-\beta}
                                      e^{u+zu^{-\alpha}}
                             du
\end{equation}
where ($Ha$) denotes Hankel's contour.

To prove (\ref{W-integral-representation}), let us write
the integrated function in the form of a power series in $z$
and perform term-by-term integration using the well-known formula
\cite[1.6(2)]{HTF-1}:
$$
        \frac{1}{\Gamma (z)} = \int_{Ha} e^{u} u^{-z} du .
$$

\noindent
\underline{{\bf Relation to other functions.}}
It follows from the definition (\ref{W-Definition}) that
\begin{equation}
      W(z; 0, 1) = e^z
\end{equation}
\begin{equation}
      \left(
           \frac{z}{2}
      \right)^\nu
      W(\mp \frac{z^2}{4}; 1, \nu + 1) = \left\{
                                \begin{array}{l}
                                J_{\nu}(z) \\
                                I_{\nu}(z)
                                \end{array}
                             \right.
\end{equation}

Taking $\beta = 1 -\alpha$, we obtain Mainardi's function $M(z; \alpha)$:
\begin{equation} \label{Mainardi-Definition}
      W(-z; -\alpha, 1-\alpha) =
               M(z; \alpha) = \sum_{k=0}^{\infty}
                              \frac{(-1)^{k} z^{k}}
                                   {k! \,\, \Gamma (-\alpha (k+1) + 1)}
\end{equation}

The following particular case of the Wright function
was given by Mainardi \cite{Mainardi-WASCOM}:
\begin{equation}
      W(-z; -\frac{1}{2}, -\frac{1}{2}) = M(z; \frac{1}{2}) =
                  \frac{1}{\sqrt{\pi}}\exp{\left(-\frac{z^2}{4}\right)}.
\end{equation}

We see that the Wright function is a generalization of
the exponential function and the Bessel functions.
For $\alpha > 0$ and  $\beta > 0$ it is an entire function in $z$
\cite{HTF-3}.

Recently Mainardi \cite{Mainardi-WASCOM}
pointed out that $W(z; \alpha, \beta)$
is an entire function in $z$ also for \mbox{$-1 <\alpha < 0$}.

Let us  prove this statement. Using the well-known relationship
\begin{equation}
      \Gamma (y) \Gamma (1-y) = \frac{\pi}{ \sin (\pi y) },
\end{equation}
we can write the Wright function in the form
\begin{equation}
      W(z; \alpha, \beta) = \frac{1}{\pi}
                            \sum_{k=0}^{\infty}
           \frac{z^{k} \Gamma (1 -\alpha k -\beta) \sin \pi (\alpha k + \beta)}
                    {k!}
\end{equation}

Let us introduce an auxilary majorizing series
\begin{equation} \label{Mazoranta}
      S = \frac{1}{\pi} \sum_{k=0}^{\infty}
             \left|
                \frac{\Gamma (1-\alpha k -\beta)}{k!}
             \right|
             |z|^k
\end{equation}
The convergence radius of series (\ref{Mazoranta}) for \mbox{$-1 <\alpha <0$}
is infinite:
\begin{equation}
      R = \lim_{k \rightarrow \infty} \left|
                            \frac{\Gamma (1- \alpha k-\beta)}{k!}
                            \frac{(k+1)!}{\Gamma (1-\alpha k - \alpha -\beta)}
                                      \right|
        = \lim_{k \rightarrow \infty} \frac{k+1}{|\alpha|^\alpha k^{-\alpha}}
        = \infty.
\end{equation}
(We use here relationship \cite[formula 1.18(4)]{HTF-1}.)

It follows the comparison of the series (\ref{W-Definition})
and (\ref{Mazoranta}) that
for $\alpha > -1$ and arbitrary $\beta$
the convergence radius of the series representation of
the Wright function $W(z; \alpha, \beta)$ is infinite,
and the Wright function is an entire function.

The Wright function plays an important role in the solution of
linear partial fractional differential equations, e.g.
fractional diffusion--wave equation.

\section{Tools for testing candidate solutions}

The following tools are necessary
for checking by substitution if candidate solutions,
obtained by any more or less formal method,
satisfy corresponding equations and initial conditions.

\vspace{1.5ex}
\noindent
\underline{{\bf Fractional derivatives of the Mittag-Leffler function.}}\\
By fractional-order differentiation $_{0}D_{t}^{\gamma}$
($\gamma$ is an arbitrary real number)
of series representation (\ref{ML-Definition}) we obtain
\begin{equation} \label{Derivative-of-E}
_{0}D_{t}^{\gamma}(t^{\alpha k +\beta -1}
                   E_{\alpha, \beta}^{(k)}(\lambda t^{\alpha})) =
      t^{\alpha k +\beta - \gamma -1}
      E_{\alpha, \beta - \gamma}^{(k)}(\lambda t^{\alpha})
\end{equation}

The particular case of
relationship (\ref{Derivative-of-E}) for $k=0$ and integer $\gamma$
is given in \cite{HTF-3}, equation 18.1(25).

\vspace{1.5ex}
\noindent
\underline{{\bf
A rule for fractional differentiation
of an integral depending on}}

\noindent
\underline{{\bf a  parameter}},when the upper limit
also depends on the parameter:
\begin{equation} \label{Parameter}
     _{0}D_{t}^{\alpha}
     \int_{0}^{t} \!\! K(t,\tau)d\tau =
     \int_{0}^{t} \!\!\, _{\tau}D_{t}^{\alpha}K(t, \tau) d\tau +
          \lim_{\tau \rightarrow t-0} \, _{\tau}D_{t}^{\alpha-1}K(t, \tau),
          \hspace{1em} (0 < \alpha < 1).
\end{equation}

The proof is straightforward:
\begin{eqnarray*}
   _{0}D_{x}^{\alpha} \int_{0}^{x} \!\! K(x,\xi)d\xi  & = &
      \frac{1}{\Gamma(1-\alpha)} \frac{d}{dx}
       \int_{0}^{x} \frac{dt}{(x-t)^{\alpha}}
        \int_{0}^{t} \!\! K(t,\xi)d\xi \\
       & = &
         \frac{1}{\Gamma(1-\alpha)} \frac{d}{dx}
         \int_{0}^{x} \!\! d\xi
         \int_{\xi}^{x} \!\! \frac{K(t, \xi) dt}{(x-t)^{\alpha}} \\
    & = & \frac{d}{dx}
                   \int_{0}^{x} \!\! \widetilde{K}(x,\xi) d\xi =
                   \int_{0}^{x} \!\! \frac{d}{dx} \widetilde{K}(x,\xi) d\xi +
                   \lim_{\xi \rightarrow x-0} \widetilde{K}(x, \xi) \\
     & = & \int_{0}^{x} \!\!\, _{\xi}D_{x}^{\alpha}K(x, \xi) d\xi +
          \lim_{\xi \rightarrow x-0} \, _{\xi}D_{x}^{\alpha-1}K(x, \xi),
\end{eqnarray*}
where
$$
   \widetilde{K}(x,\xi) = \frac{1}{\Gamma(1-\alpha)}
         \int_{\xi}^{x} \!\! \frac{K(t, \xi) dt}{(x-t)^{\alpha}}.
$$

The following important particular case must be mentioned.
If we have \mbox{$K(t -\tau)f(\xi)$} instead of $K(t, \tau)$,
then relationship (\ref{Parameter}) takes the form:
\begin{equation}\label{Convolution-Kernel}
 _{0}D_{t}^{\alpha} \int_{0}^{t} \!\! K(t -\tau) f(\tau) d\tau =
    \int_{0}^{t} \!\! _{0}D_{\tau}^{\alpha}K(\tau)f(t-\tau)d\tau +
     \lim_{\tau \rightarrow +0} f(t-\tau) \, _{0}D_{\tau}^{\alpha -1}K(\tau))
\end{equation}

\chapter{Standard fractional differential equations}

The following examples illustrate the use of (\ref{Laplace-Transform})
for solving fractional-order differential equations
with constant coefficients.
In this chapter we use the classic formula for the Laplace transform
of the fractional derivative, as given, e.g., in \cite{OS}, p.134
or \cite{Miller-Ross}, p.123:
\begin{equation}\label{Classics}
      \int_{0}^{\infty} \,\, e^{-pt} \,_{0}D_{t}^{\alpha}f(t) \,\, dt =
      p^{\alpha}F(p) - \sum_{k=0}^{n-1} p^{k} \,\,
                    \left.  _{0}D_{t}^{\alpha - k - 1}f(t) \right|_{t=0},
      \hspace{1em}
      (n-1 < \alpha \leq n).
\end{equation}

\section{Ordinary linear fractional \protect \\ differential equations}

In this section we give some examples of solution of ordinary
linear differential equations of the fractional order.

%
%
\example \label{Oldham-Spanier-1}
A slight generalization of an equation solved in (\cite{OS}, p.157):
\begin{equation}\label{Example-1}
        _{0}D_{t}^{1/2}f(t) + af(t) = 0, \hspace{1em} (t>0);
        \hspace{2em}
        \left. _{0}D_{t}^{-1/2}f(t) \right|_{t=0}=C
\end{equation}

Applying the Laplace transform we obtain
$$
F(p) = \frac{C}{p^{1/2} + a},
         \hspace{3em} \left. C= \, _{0}D_{t}^{-1/2}f(t) \right|_{t=0}
$$
and the inverse transform with a help of (\ref{Laplace-Transform-Semi})
gives the solution of (\ref{Example-1}):
\begin{equation}\label{Solution-1}
f(t) = C t^{-1/2} E_{\frac{1}{2},\frac{1}{2} }(-a\sqrt{t}).
\end{equation}
Using series expansion (\ref{ML-Definition}) of $E_{\alpha, \beta}(t)$,
it is easy to check that for $a=1$ solution (\ref{Solution-1})
is identical to the solution
$f(t)=C(\frac{1}{\sqrt{\pi t}} - e^{t} \mbox{erfc}(\sqrt{t}))$,
obtained in \cite{OS} in a more complex way.

\example \label{Oldham-Spanier-2}
Let us consider the following equation:
\begin{equation}\label{Example-2}
        _{0}D_{t}^{Q}f(t) + \, _{0}D_{t}^{q}f(t) = h(t),
\end{equation}
which ''encounters very great difficulties except when the difference
$q-Q$ is integer or half-integer'' (\cite{OS}, p.156).

Suppose that $0 < q < Q < 1$.
Laplace transform of (\ref{Example-2}) leads to
\begin{equation}
(p^Q + p^q)F(p) = C + H(p), \hspace{1em}
      \mbox{where} \hspace{1em} C = \left.
                                     \left(
                                      _{0}D_{t}^{q-1}f(t) + \, _{0}D_{t}^{Q-1}
                                     \right)
                                     \right|_{t=0},
\end{equation}
and
\begin{equation}
F(p) = \frac{C+H(p)}{p^Q + p^q} = \frac{C+H(p)}{p^q (p^{Q-q}+1)} =
     \left(
        C+H(p)
     \right)
        \frac{p^{-q}}{p^{Q-q}+1}.
\end{equation}

After inversion with a help of (\ref{Laplace-Transform})
for $\alpha=Q-q$ and $\beta=Q$,
we obtain the solution:

\begin{equation}
f(t) = C \, G(t)
         +
         \int_{0}^{t}
         G(t-\tau)
         h(\tau) d\tau,
\end{equation}
$$
    C = \left.
              \left(
                  _{0}D_{t}^{q-1}f(t) + \, _{0}D_{t}^{Q-1}f(t)
              \right)
        \right|_{t=0},
                        \hspace{1em}
     G(t) = t^{Q-1}E_{Q-q, Q}(-t^{Q-q})
$$

The case $0< q< Q < n$ (for example, equation obtained in \cite{Oustaloup})
can be solved similarly.

\example \label{Samko}
Let us consider the following initial value problem
for a non-homogeneous fractional differential equation
under non-zero initial conditions:
\begin{equation}
   _{0}D_{t}^{\alpha} y(t) - \lambda y(t) = h(t),
       \hspace{1em}  (t>0);   \label{Example-4}
\hspace{1em}  \left.
        _{0}D_{t}^{\alpha -k} y(t)
      \right|_{t=0} = b_{k},
      \hspace{1em}  (k = 1, 2, \ldots, n) \label{Init-Cond-4}
\end{equation}
where $n-1 < \alpha < n$.
Problem (\ref{Example-4}) was analytically solved in \cite{SKM}
by the iteration method.
With the help of the Laplace transform and
formula (\ref{Laplace-Transform})
we obtain the same solution directly and more quickly.

Indeed, taking into account initial conditions (\ref{Init-Cond-4}),
the Laplace transform of equation (\ref{Example-4}) yields
$$
   p^{\alpha} Y(p) - \lambda Y(p) = H(p) + \sum_{k=1}^{n} b_{k}p^{k-1},
$$
\begin{equation}
   Y(p) = \frac{\mbox{$H(p)$}}{\mbox{$p^{\alpha} - \lambda$}} +
         \sum_{k=1}^{n} b_{k} \frac{p^{k-1}}{p^{\alpha} - \lambda }
                                 \label{Transform-4}
\end{equation}
and the inverse transform using (\ref{Laplace-Transform}) gives
the solution:
\begin{equation} \label{Solution-4}
y(t) = \sum_{k=1}^{n} b_{k} t^{\alpha - k}
                       E_{\alpha, \alpha - k + 1} (\lambda t^{\alpha}) +
       \int_{0}^{t} (t - \tau)^{\alpha - 1}
                  E_{\alpha, \alpha } (\lambda (t - \tau)^{\alpha})
                  h(\tau) d\tau.
\end{equation}

\newpage

\section{Partial linear fractional \protect \\ differential equations}

The proposed approach can be successfully used for
solving partial linear differential equations of the fractional order.

\example \label{Nigmatullin}
{\bf Nigmatullin's fractional diffusion equation.}

Let us consider the following initial boundary value problem for the
fractional diffusion equation in one space dimension:
\begin{equation}
   _{0}D_{t}^{\alpha}u(x,t) =
             \lambda^{2} \frac{\partial^{2}u(x,t)}{\partial x^{2}},
                        \hspace{1em}
                       ( t>0, \hspace{0.5em} -\infty < x < \infty )
                                                         \label{Example-3};
\end{equation}
\begin{equation}
   \hspace{1em}
   \lim_{x \rightarrow \pm \infty} u(x,t) = 0; \label{Boundary-Conditions}
   \hspace{1em}
   \left. _{0}D_{t}^{\alpha -1}u(x,t) \right|_{t=0}=\varphi (x).
                \label{Initial-Condition}
\end{equation}
We assume here $0 < \alpha < 1$.
Equation of the type (\ref{Example-3}) was deduced
by Nigmatullin \cite{Nigmatullin-Transfer}
and by Westerlund \cite{Westerlund}
and studied by Mainardi \cite{Mainardi-WASCOM}).
We will give a simple solution of problem (\ref{Example-3})
demonstrating once again the advantage of using
the Mittag-Leffler function in two parameters (\ref{ML-Definition}).

Taking into account boundary conditions (\ref{Boundary-Conditions}),
the Fourier transform with respect to variable $x$ gives:
\begin{eqnarray}
    & &
        _{0}D_{t}^{\alpha} \overline{u}(\beta, t) +
                  \lambda^{2} \beta^{2} \overline{u}(\beta,t) = 0
                                            \label{F-1}  \\
    & &
   \left. _{0}D_{t}^{\alpha -1}\overline{u}(x,t) \right|_{t=0}=
                                         \overline{\varphi} (\beta),
                                            \label{F-2}
\end{eqnarray}
where $\beta$ is the Fourier transform parameter.
Applying Laplace transform to (\ref{F-1})
and using initial condition (\ref{F-2}) we obtain
\vspace{-1ex}
\begin{equation} \label{Transformanta}
        \overline{U}(\beta, p) =
              \frac{\varphi (\beta)}{p^\alpha + \lambda^{2}\beta^{2}}.
\end{equation}

Inverse Laplace transform of (\ref{Transformanta})
using (\ref{Laplace-Transform}) gives
\begin{equation}
  \overline{u}(\beta,t) = \varphi (\beta) t^{\alpha-1}
               E_{\alpha, \alpha} (-\lambda^{2} \beta^{2} t^\alpha),
\end{equation}
and inverse Fourier transform produces the solution
of the problem (\ref{Example-3})-(\ref{Initial-Condition}):

\begin{equation} \label{Solution-3}
u(x,t)  =  \int_{-\infty}^{\infty} G(x-\xi, t) \varphi (\xi) d\xi,
\hspace{3em}
\label{Kernel-3}
\end{equation}
\begin{equation}\label{Nigmatullin-G}
G(x,t) = \frac{1}{\pi}
         \int_{0}^{\infty}
         t^{\alpha-1} E_{\alpha, \alpha} (-\lambda^{2} \beta^{2} t^\alpha)
         \cos{\beta x}
         d\beta.
\end{equation}

Let us evaluate integral (\ref{Nigmatullin-G}).
The Laplace transform of (\ref{Nigmatullin-G})
and formula 1.2(11) from \cite{TIT-I} produce
\begin{equation}
   g(x,p) = \frac{1}{\pi}
              \int_{0}^{\infty}
               \frac{\cos(\beta x \,\, d\beta)}
                    {\lambda^{2}\beta^{2} + p^\alpha}
          = \frac{1}{2\lambda}p^{-\alpha /2}e^{-|x|\lambda^{-1}p^{\alpha /2}},
\end{equation}
and the Laplace inversion gives
\begin{equation}
      G(x,t) = \frac{1}{4 \lambda \pi i} \int \limits_{Br} e^{pt}
                p^{-\frac{\alpha}{2}} \exp (-|x|\lambda^{-1} p^{\alpha /2})
                dp.
\end{equation}

Performing substitutions $\sigma = pt$ and $z= |x|\lambda^{-1}t^{-\rho}$
($\rho = \alpha /2$) and transformation of
the Bromwich contour ($Br$) to the Hankel contour ($Ha$),
as was done in a similar case by Mainardi \cite{Mainardi-WASCOM},
we obtain
\begin{equation} \label{Nigmatullin-G-Wright}
      G(x,t) = \frac{t^{1-\rho}}{2\lambda}
               \frac{1}{2 \pi i}
               \int \limits_{Ha}
               e^{\sigma - z \sigma^\rho } \frac{d\sigma}{\sigma^{\rho}}
             = \frac{1}{2\lambda}t^{\rho -1} W(-z, -\rho, \rho),
             \qquad
             z= \frac{|x|}{\lambda t^\rho}
\end{equation}
where $W(z, \lambda, \mu)$ is the Wright function (\ref{W-Definition}).
We would like to note that, in fact, we  have just evaluated
the Fourier cosine-transform of the function
$u_1(\beta) = t^{\alpha-1} E_{\alpha, \alpha} (-\lambda^{2} \beta^{2} t^\alpha)$.

It is easy to check that for $\alpha = 1$ (traditional diffusion equation)
fractional Green function (\ref{Nigmatullin-G-Wright})
reduces to the classic expression
\begin{equation}\label{Classic}
G(x,t) = \frac{1}{2 \lambda \sqrt{\pi t}}
         \exp (- \frac{x^2}{4 \lambda^{2} t} ).
\end{equation}

\example \label{Wyss}
{\bf Wyss's fractional diffusion equation.}

The following example shows the effectiveness of the proposed method also
for fractional integral equations. Let us consider the Wyss's type formulation
of the diffusion equation \cite{Schneider-Wyss}
(for simplicity and comparison with the previous
example --- in one space dimension):
\begin{equation}
u(x, t) = \varphi(x) + \lambda^2 \, _{0}D_{t}^{-\alpha}
                \frac{\partial^2 u(x,t)}{\partial x^2},
                \hspace{1em} ( -\infty < x < \infty, \hspace{1em} t>0);
                \label{Example-5}
\end{equation}
\begin{equation}
\lim_{x \rightarrow \pm \infty} u(x,t) = 0,
                \label{Boundary-Condition-5}
\hspace{1em}
u(x, 0) =\varphi(x),
                \label{Initial-Condition-5}
\end{equation}

Applying the Fourier transform with respect to space variable $x$
and the Laplace transform with respect to time $t$, we obtain:
\begin{equation}
\overline{U}(\beta, p) = \frac{\varphi(\beta) p^{\alpha - 1}}
                              {p^\alpha + \lambda^2 \beta^2},
\end{equation}
where $\overline{U}(\beta, p)$ is the Fourier-Laplace transform of $u(x,t)$,
$\beta$ is the Fourier transfrom parameter
and $p$ is the Laplace transform parameter.

Inverting Laplace and Fourier transforms as it was done
in Example \ref{Nigmatullin}, we obtain the solution
of problem (\ref{Example-5}):
\begin{equation}\label{Solution-5}
u(x, t) = \int_{-\infty}^{\infty} G(x-\xi, t) \varphi (\xi) d\xi,
\end{equation}
\begin{equation}\label{Wyss-G}
G(x,t) = \frac{1}{\pi}
         \int_{0}^{\infty}
         E_{\alpha, 1} (-\lambda^{2} \beta^{2} t^\alpha)
         \cos{\beta x}
         d\beta.
\end{equation}

Let us evaluate integral (\ref{Wyss-G}).
The Laplace transform of (\ref{Wyss-G})
and formula 1.2(11) from \cite{TIT-I} produce
\begin{equation}
      g(x, p) =  \frac{p^{\alpha -1}}{\pi} \int \limits_{0}^{\infty}
                \frac{ \cos (\beta x) d\beta}{p^\alpha + \lambda^{2}\beta^2}
             = \frac{1}{2\lambda}
               p^{\frac{\alpha}{2} -1} \exp (-|x| \lambda^{-1} p^{\alpha /2}),
\end{equation}
and the the inverse Laplace transform gives:
\begin{equation}
      G(x,t) = \frac{1}{4 \lambda \pi i} \int \limits_{Br} e^{pt}
                p^{\frac{\alpha}{2} -1} \exp (-|x|\lambda^{-1} p^{\alpha /2})
                dp.
\end{equation}

Performing substitutions $\sigma = pt$ and $z= |x|\lambda^{-1}t^{-\rho}$
($\rho = \alpha /2$) and transformation of
the Bromwich contour ($Br$) to the Hankel contour ($Ha$),
as it was done in a similar case by Mainardi \cite{Mainardi-WASCOM},
we obtain
\begin{equation} \label{Wyss-G-Wright}
      G(x,t) = \frac{t^{-\rho}}{2\lambda}
               \frac{1}{2 \pi i}
               \int \limits_{Ha}
               e^{\sigma - z \sigma^\rho } \frac{d\sigma}{\sigma^{1-\rho}}
             = \frac{1}{2\lambda}t^{-\rho} M(z, \rho),
             \quad
             z= \frac{|x|}{\lambda t^\rho}
\end{equation}
where $M(z, \rho)=W(-z, -\rho, 1-\rho)$ is the Mainardi function
(\ref{Mainardi-Definition}).

The last expression is identical to the expression, obtained by Mainardi
\cite{Mainardi-WASCOM} in another way.

We would like to note at this point, as in the previous example,
that we have just evaluated the Fourier cosine-transform of the function
$u_2(\beta) = E_{\alpha, 1} (-\lambda^{2} \beta^{2} t^\alpha)$.

For $\alpha = 1$ the fractional Green's function (\ref{Wyss-G-Wright})
also reduces
to the classic expression (\ref{Classic}). The case of arbitrary
number of space dimensions can be solved similarly.

For $\alpha = 1$ both generalizations (Nigmatullin's as well as Wyss's)
of the diffusion problem give the standard diffusion problem,
and the solutions reduce to the classic solution.
However, it is obvious that the asymptotic behavior of (\ref{Solution-3})
and (\ref{Solution-5}) for $t \rightarrow 0$,
and $t\rightarrow \infty$ is different
(see also discussion in \cite{GN}
on two different generalizations of the standard relaxation equation
and discussion in \cite{Friedrich}
on two fractional models --- one based on fractional derivatives
and the other based on fractional integrals --- for mechanical stress relaxation).

This difference was caused by initial conditions
of different types. The class of solutions is determined by
the number and the type of initial conditions.

\chapter{Sequential fractional differential equations}

\section{The Laplace transform of a sequential \protect \\
                fractional differential operator}

Let us consider the following initial value problem:
\begin{equation} \label{IVP}
    _{0}{\cal L}_{t}y(t)= f(t);
    \hspace{2em}
    \left.
    _{0}{\cal D}_{t}^{\sigma_{k} -1}y(t)
    \right|_{t=0} = b_{k},
    \hspace{1em}
    k= 1, \, \ldots \, , n
\end{equation}
\begin{equation} \label{LFDO}
     _{a}{\cal L}_{t} y(t) \equiv \,\,
    _{a}{\cal D}_{t}^{\sigma _{n}}y(t) +
    \sum_{k=1}^{n-1}p_{k}(t) \, _{a}{\cal D}_{t}^{\sigma _{n-k}}y(t) +
    p_{n}(t)y(t),
\end{equation}
$$
 _{a}{\cal D}_{t}^{\sigma_{k}} \equiv  \, _{a}D_{t}^{\alpha_{k}}
                                       \, _{a}D_{t}^{\alpha_{k-1}} \ldots
                                       \,_{a}D_{t}^{\alpha_{1}};
\hspace{1em}
 _{a}{\cal D}_{t}^{\sigma_{k} -1} \equiv  \, _{a}D_{t}^{\alpha_{k}-1}
                                       \, _{a}D_{t}^{\alpha_{k-1}} \ldots
                                       \,_{a}D_{t}^{\alpha_{1}};
$$
$$
\sigma _{k} =\sum_{j=1}^{k}\alpha _{j}, \hspace{1em} (k= 1, 2, \ldots, n);
\hspace{1.5em} 0 \leq \alpha _{j} \leq 1, \hspace{1em} (j = 1, 2, \ldots,n).
$$

Fractional differential equation (\ref{IVP})
is a sequential fractional
differential equation,
according to the terminology used by Miller and Ross \cite{Miller-Ross}.
To extend the Laplace transform method
using the advantages of (\ref{Laplace-Transform}) for such equations
with constant coefficients,
we obtained the following formula:

\vspace{1.5ex}
\noindent
\underline{{\bf Laplace transform of the
sequential fractional derivative $_{0}{\cal D}_{t}^{\sigma_{m}}y(t)$.}}
\begin{equation}\label{FLT}
\int_{0}^{\infty} e^{-pt} \, _{0}{\cal D}_{t}^{\sigma_{m}}f(t) dt =
p^{\sigma_{m}}F(p) -
\sum_{k=0}^{m-1}p^{\sigma_{m}-\sigma_{m-k}}
            \left.
            \, _{0}{\cal D}_{t}^{\sigma_{m-k}-1}f(t)
            \right|_{t=0},
\end{equation}
\begin{displaymath}
_{a}{\cal D}_{t}^{\sigma_{m - k}-1}
                               \equiv  \, _{a}D_{t}^{\alpha_{m-k}-1}
                                       \, _{a}D_{t}^{\alpha_{m-k-1}} \ldots
                                       \, _{a}D_{t}^{\alpha_{1}}.
\end{displaymath}

The particular case of (\ref{FLT}) for $f(t)$ m-times differentiable,
$\alpha_{m}=\mu$, $\alpha_{k}=1$, $(k=1,2,\ldots,m-1)$ was obtained
by Caputo \cite[p.41]{Caputo-ED} much earlier.
Taking $\alpha_{1}=\mu$, $\alpha_{k}=1$, $(k=2,3, \ldots, m)$
leads under obvious assumptions to the classic formula
(\ref{Classics}).

\section{Ordinary linear fractional \protect \\ differential equations}

In this section we give solutions of the ''sequential'' analogues
of the ''standard'' linear ordinary fractional differential equations
with constant coefficients.
Of course, we must take appropriate initial conditions, also
in terms of sequential fractional derivatives.

\example \label{Example-1a}
Let us consider the sequential analogue of Example \ref{Oldham-Spanier-1}:
\begin{equation} \label{Equation-1a}
      _{0}D_{t}^{\alpha} \left( \,\, _{0}D_{t}^{\beta} y(t)
                         \right)
                         + a y(t) = 0
\end{equation}
\begin{equation} \label{Initial-Conditions-1a}
      \left. _{0}D_{t}^{\alpha -1} \left( \,\, _{0}D_{t}^{\beta} y(t)
                         \right) \right|_{t=0} = b_1,
      \qquad
      \left. _{0}D_{t}^{\beta -1} y(t) \right|_{t=0} = b_2,
\end{equation}
where $0 < \alpha < 1$, $0< \beta <1$, $\alpha +\beta = 1/2$.

Our formula (\ref{FLT}) of the Laplace transform of the sequential
fractional derivative allows us to utilize
the initial conditions (\ref{Initial-Conditions-1a}).
To use (\ref{FLT}), we take $\alpha_1 = \alpha$, $\alpha_2 = \beta$
and $m=2$. Therefore, $\sigma_1 =\alpha$, $\sigma_2 =\alpha + \beta$.
Then the Laplace transform (\ref{FLT})  of equation (\ref{Equation-1a}) gives:
\begin{equation}
      (p^{\alpha +\beta} + a) Y(p) = p^\beta b_2 + b_1,
\end{equation}
\begin{equation}
      Y(p) = b_2 \frac{p^{\beta}}{p^{\alpha +\beta}+a} +
             b_1 \frac{1}{p^{\alpha +\beta}+a},
\end{equation}
and after the Laplace inversion with the help
of (\ref{Laplace-Transform}) we find the solution to the problem
(\ref{Equation-1a})-(\ref{Initial-Conditions-1a}):
\begin{equation} \label{Solution-1a}
   y(t) = b_2 t^{\alpha-1} E_{\alpha +\beta, \alpha}(-at^{\alpha +\beta}) +
          b_1 t^{\alpha +\beta -1}E_{\alpha +\beta, \alpha +\beta}
                                    (-at^{\alpha +\beta}).
\end{equation}
For $\beta = 0$ and $\alpha = 1/2$ (and assuming, of course, $b_2 =0$),
we can obtain from (\ref{Solution-1a})
the solution of Example \ref{Oldham-Spanier-1}.

\example \label{Oldham-Spanier-2-seq}
Let us now consider the following sequential analogue
for the equation, considered in Example \ref{Oldham-Spanier-2}:
\begin{equation} \label{Equation-2a}
        _{0}D_{t}^{\alpha} \left(
                                \,\, _{0}D_{t}^{\beta}y(t)
                           \right)
               + \,\, _{0}D_{t}^{q}y(t) =h(t),
\end{equation}
where $0 < \alpha < 1$, $0< \beta <1$, $\alpha +\beta = Q > q$.

The Laplace transform (\ref{FLT}) of equation (\ref{Equation-2a}) gives:
\begin{equation}
       (p^{\alpha +\beta} + p^q) Y(p) = H(p) + p^\beta b_2 + b_1,
\end{equation}
$$
   b_1 = \,\, _{0}D_{t}^{\alpha -1} \left.
                                       \left(
                                          \,\, _{0}D_{t}^{\beta}y(t)
                                       \right)
                                   \right|_{t=0} +
         \left. \,\, _{0}D_{t}^{q-1}y(t) \right|_{t=0},
$$
$$
   b_2 = \left. \,\, _{0}D_{t}^{\beta -1}y(t) \right|_{t=0}.
$$

Writing $Y(p)$ in the form
\begin{equation}
      Y(p) = \frac{p^{-q}H(p)}{p^{\alpha +\beta - q}+1} +
             b_2 \frac{p^{\beta -q}}{p^{\alpha +\beta - q}+1} +
             b_1 \frac{p^{-q}}{p^{\alpha +\beta - q}+1}
\end{equation}
and finding the inverse Laplace transform with the help
of (\ref{Laplace-Transform}), we obtain the solution:
\begin{eqnarray}
    y(t) & = & b_2 t^{\alpha -1} E_{\alpha +\beta -q, \alpha}
                                   (-t^{\alpha +\beta -q}) +
               b_1 t^{\alpha +\beta -q} E_{\alpha +\beta -q, \alpha +\beta}
                                   (-t^{\alpha +\beta -q}) +
                                   \nonumber \\
        & &      + \int_{0}^{t}
                       (t-\tau)^{\alpha +\beta -1}
                       E_{\alpha +\beta -q, \alpha +\beta}
                          \left(
                                -(t-\tau)^{\alpha +\beta -q}
                          \right)
               h(\tau) d\tau.   \label{Solution-2a}
\end{eqnarray}

It is easy to see that this solution contains
solution of Example \ref{Oldham-Spanier-2} as a special case.

\example \label{Example-3a}
Let us consider the following initial value problem
for the sequential fractional differential equation:

\begin{equation}\label{Example-6}
       \hspace{-0.5em}
        _{0}D_{t}^{\alpha_{2}}(\, _{0}D_{t}^{\alpha_{1}}y(t))
                   - \lambda y(t)= h(t);
        \hspace{1em}
        (0< \alpha_1, \alpha_2 < 1)
\end{equation}
\begin{equation}
        \left. \, _{0}D_{t}^{\alpha_{2}-1}
                  \, _{0}D_{t}^{\alpha_{1}} y(t)\right|_{t=0} = b_{1},
        \hspace{0.5em}
        \left. \, _{0}D_{t}^{\alpha_{1}-1} y(t)\right|_{t=0} = b_{2};
\end{equation}

The Laplace transform (\ref{FLT}) of equation (\ref{Example-6}) gives
$$
(p^{\alpha_{1}+\alpha_{2}} - \lambda)Y(p) = p^{\alpha_{2}}b_{2} + b_{1},
$$
and after inversion using (\ref{Laplace-Transform}) we obtain
the solution:

$$
y(t) = b_{2}t^{\alpha_{1}-1}
       E_{\alpha, \alpha_{1}}(\lambda t^{\alpha}) +
       b_{1}t^{\alpha-1}
       E_{\alpha, \alpha}
       (\lambda t^{\alpha}) +
$$
\begin{equation} \label{Solution-6}
       +
       \int_{0}^{t}
       (t-\tau)^{\alpha-1}
       E_{\alpha, \alpha}
       (\lambda (t-\tau)^{\alpha})
       h(\tau)d\tau,
       \hspace{2em}
       (\alpha=\alpha_{1}+\alpha_{2})
\end{equation}

Let us take $\alpha$ the same as in Example \ref{Samko}.
Using (\ref{ML-Definition}), (\ref{Derivative-of-E})
and (\ref{Convolution-Kernel}), it is easy to
verify that (\ref{Solution-6}) is the solution of (\ref{Example-6}).
It is also worthwhile to note that
if $b_1 \neq 0$, $b_2 \neq 0$, then
(\ref{Solution-6}) is {\em not}\/ a solution of equation
$\, _{0}D_{t}^{\alpha}y(t) - \lambda y(t)=h(t)$ from Example \ref{Samko};
also (\ref{Solution-4}) is {\em not}\/
a solution of equation (\ref{Example-6}).
On the other hand, equations (\ref{Example-4}) and (\ref{Example-6})
are very close one to another: the fractional Green's function
in both cases is
$G(t)=(t)^{\alpha-1}E_{\alpha, \alpha}(\lambda t^{\alpha})$.

\section{Partial linear fractional \protect \\ differential equations}

\example \label{Mainardi-problem}
Let us consider Mainardi's \cite{Mainardi-WASCOM}
initial value problem for the fractional diffusion--wave equation:
\begin{equation}\label{M-equation}
        _{0}D_{t}^{\alpha}u(x,t) = \lambda^{2}
                      \frac{\partial^{2}u(x,t)}{\partial x^2},
        \hspace{3em}
        ( |x|<\infty, t>0 )
\end{equation}
\begin{equation}\label{M-initial-condition}
      u(x,0) = f(x),
        \hspace{3em}
        ( |x|<\infty )
\end{equation}
\begin{equation}  \label{M-boundary-condition}
      \lim_{x \pm \infty} u(x,t) = 0,
        \hspace{3em}
        ( t>0 )
\end{equation}
where $0< \alpha <1$.

The Laplace transform of problem
(\ref{M-equation})--(\ref{M-boundary-condition})
using (\ref{FLT}) for $k=2$, $\alpha_{1}=\alpha -1$, $\alpha_{2}=1$
(this gives Caputo's formula \cite{Caputo-ED}), i.e.,
\begin{equation}\label{Laplace-Caputo}
    {\cal L}\{\,\, _{0}D_{t}^{\alpha} y(t)\} = p^\alpha Y(p) - p^{\alpha -1}y(0)
\end{equation}
produces:
\begin{equation} \label{M-L-equation}
      p^{\alpha} \overline{u} (x,p) - p^{\alpha -1}f(x) =
        \lambda^{2} \frac{\partial^{2} \overline{u}(x,p)}{\partial x^2},
      \hspace{3em}
      (|x| < \infty)
\end{equation}
\begin{equation} \label{M-L-boundary-conditions}
      \lim_{x \pm \infty} \overline{u}(x,p) = 0,
        \hspace{3em}
        ( t>0 )
\end{equation}

Applying now the Fourier exponential transform
to equation (\ref{M-L-equation}) and utilizing
boundary conditions (\ref{M-L-boundary-conditions}),
we obtain:
\begin{equation}
  U(\beta, p) = \frac{p^{\alpha -1}}{p^\alpha + \lambda^{2}\beta^2}F(\beta),
\end{equation}
where $U(\beta , p)$ and $F(\beta)$ are the Fourier transforms
of $\overline{u}(x,p)$  and $f(x)$.

The inverse Laplace transform of the fraction
$\frac{p^{\alpha -1}}{p^\alpha + \lambda^{2}\beta^2}$
is $E_{\alpha, 1}(- \lambda^{2}\beta^2 t^\alpha)$
($E_{\lambda, \mu}(z)$ is the Mittag-Leffler function
in two parameters). Therefore, the inversion of the Fourier
and the Laplace transform gives the solution in the following form:
\begin{equation} \label{M-solution}
      u(x,t) = \int \limits_{-\infty}^{\infty} G(x-\xi, t)f(\xi) d\xi,
\end{equation}
\begin{equation} \label{M-Green-function}
      G(x,t) = \frac{1}{\pi} \int \limits_{0}^{\infty}
           E_{\alpha, 1}(-\lambda^{2}\beta^2 t^\alpha) \cos(\beta x) d\beta
            = \frac{1}{2\lambda}t^{-\rho} W(-z, -\rho, 1-\rho),
\end{equation}
where $W(z, \lambda, \mu)$ is the Wright function (\ref{W-Definition}).
This solution is identical to the solution of the Wyss
fractional (integral) diffusion equation (\ref{Solution-5}).

\chapter{Fractional Green's function}

\section{Definition and some properties}

In this section we consider equation (\ref{IVP})
under homogeneous initial conditions $b_{k}=0$, $(k=1,\ldots,n)$, i.e.
\begin{equation} \label{IVP-G}
    _{0}{\cal L}_{t}y(t)= f(t);
    \hspace{2em}
    \left.
    _{0}{\cal D}_{t}^{\sigma_{k} -1}y(t)
    \right|_{t=0} = 0,
    \hspace{1em}
    k= 1, \, \ldots \, , n
\end{equation}
$$
     _{a}{\cal L}_{t} y(t) \equiv \,\,
    _{a}{\cal D}_{t}^{\sigma _{n}}y(t) +
    \sum_{k=1}^{n-1}p_{k}(t) \, _{a}{\cal D}_{t}^{\sigma _{n-k}}y(t) +
    p_{n}(t)y(t),
$$
$$
 _{a}{\cal D}_{t}^{\sigma_{k}} \equiv  \, _{a}D_{t}^{\alpha_{k}}
                                       \, _{a}D_{t}^{\alpha_{k-1}} \ldots
                                       \,_{a}D_{t}^{\alpha_{1}};
\hspace{1em}
 _{a}{\cal D}_{t}^{\sigma_{k} -1} \equiv  \, _{a}D_{t}^{\alpha_{k}-1}
                                       \, _{a}D_{t}^{\alpha_{k-1}} \ldots
                                       \,_{a}D_{t}^{\alpha_{1}};
$$
$$
\sigma _{k} =\sum_{j=1}^{k}\alpha _{j}, \hspace{1em} (k= 1, 2, \ldots, n);
\hspace{1.5em} 0 \leq \alpha _{j} \leq 1, \hspace{1em} (j = 1, 2, \ldots,n).
$$

The following definition is a ''fractionalization'' of the definition
given in \cite{Naimark}.

\vspace{1.5ex}
\noindent
\underline{{\bf Definition.}}
Function $G(t, \tau)$ satisfying  the following conditions

a) $_{\tau}{\cal L}_{t}G(t, \tau)=0$
    for every $\tau \in (0, t)$;

b) $\lim_{\tau \rightarrow t-0}
      (\, _{\tau}{\cal D}_{t}^{\sigma _{k}-1} G(t, \tau)) = \delta_{k,n}$,
       \hspace{1em}
       $k=0, 1, \ldots, n$, \\
       \hspace*{2em}
       ($\delta_{k,n}$\/ is Kronecker's delta);

c) $\lim_{
          \mbox{$
          \begin{array}{ll}
               \mbox{{\scriptsize $ \tau, t \rightarrow +0 $}} \\
               \mbox{{\scriptsize $ \tau < t $}}
          \end{array}
          $}
          }
 ( \,_{\tau}{\cal D}_{t}^{\sigma_{k}}G(t, \tau) ) = 0,
      \hspace{1em}
      k=0, 1, \ldots, n-1$\\
is called Green's function of equation (\ref{IVP-G}).

\vspace{1.5ex}

\vspace{1.5ex}
\noindent
\underline{{\bf Properties.}}

\underline{{\bf 1.}}
Using (\ref{Parameter}), it can be shown that
$y(t) = \int_{0}^{t}G(t, \tau)f(\tau) d\tau$
is the solution of problem (\ref{IVP-G}).

\vspace{1.5ex}

Let us outline the proof of this statement.
Evaluating $_{0}D_{t}^{\sigma_1}y(t)$, $_{0}D_{t}^{\sigma_2}y(t)$,
\ldots , $_{0}D_{t}^{\sigma_n}y(t)$
using the rule (\ref{Parameter}) and condition (b) of the definition
of Green's function, we obtain:
\begin{eqnarray}
      _{0}D_{t}^{\sigma_1}y(t)
        & = & \,\, _{0}D_{t}^{\alpha_1} \int_{0}^{t} G(t, \tau) f(\tau) d\tau
                           \nonumber \\
        & = & \int_{0}^{t}\, _{\tau}D_{t}^{\alpha_1}G(t,\tau) f(\tau) d\tau
           + \lim_{\tau \rightarrow t-0}
                \,\, _{\tau}D_{t}^{\alpha_1 -1} G(t,\tau) f(\tau)
                           \nonumber \\
        & = & \int_{0}^{t}\, _{\tau}D_{t}^{\sigma_1}G(t,\tau) f(\tau) d\tau
\end{eqnarray}

\begin{eqnarray}
      _{0}D_{t}^{\sigma_2}y(t)
        & = &  \, _{0}D_{t}^{\alpha_2} (\,\, _{0}D_{t}^{\alpha_1}) =
               \,\, _{0}D_{t}^{\alpha_2}
               \int_{0}^{t} \, _{\tau}D_{t}^{\alpha_1} G(t, \tau) f(\tau) d\tau
                           \nonumber \\
        & = & \int_{0}^{t}
               \, _{\tau}D_{t}^{\alpha_2}
               (\, _{\tau}D_{t}^{\alpha_1}G(t,\tau)) f(\tau) d\tau
                           \nonumber \\
        & & \hspace*{5em}   + \lim_{\tau \rightarrow t-0}
               \, _{\tau}D_{t}^{\alpha_2 -1}
                (\, _{\tau}D_{t}^{\alpha_1} G(t,\tau)) f(\tau)
                           \nonumber \\
        & = & \int_{0}^{t}\, _{\tau}D_{t}^{\sigma_2}G(t,\tau) f(\tau) d\tau
\end{eqnarray}

\begin{eqnarray}
 \cdots & \cdots & \cdots \nonumber
\end{eqnarray}

\begin{eqnarray}
_{0}D_{t}^{\sigma_{n-1}}y(t)
      & = & \, _{0}D_{t}^{\alpha_{n-1}}(\, _{0}D_{t}^{\sigma_{n-2}}) =
             \,\, _{0}D_{t}^{\alpha_{n-1}}
             \int_{0}^{t} \, _{\tau}D_{t}^{\sigma _{n-2}} G(t, \tau) f(\tau) d\tau
             \nonumber \\
      & = & \int_{0}^{t}
               \, _{\tau}D_{t}^{\alpha_{n-1}}
               (\, _{\tau}D_{t}^{\sigma_{n-2}}G(t,\tau)) f(\tau) d\tau
                           \nonumber \\
        & & \hspace*{5em}
           + \lim_{\tau \rightarrow t-0}
               \, _{\tau}D_{t}^{\alpha_{n-1} -1}
                (\, _{\tau}D_{t}^{\sigma_{n-2}}  G(t,\tau)) f(\tau)
                           \nonumber \\
        & = & \int_{0}^{t}\, _{\tau}D_{t}^{\sigma _{n-1}}G(t,\tau) f(\tau) d\tau
\end{eqnarray}

\begin{eqnarray}
_{0}D_{t}^{\sigma_{n}}y(t)
      & = & \, _{0}D_{t}^{\alpha_{n}}(\, _{0}D_{t}^{\sigma_{n-1}}) =
             \,\, _{0}D_{t}^{\alpha_{n}}
             \int_{0}^{t} \, _{\tau}D_{t}^{\sigma _{n-1}} G(t, \tau) f(\tau) d\tau
             \nonumber \\
      & = & \int_{0}^{t}
               \, _{\tau}D_{t}^{\alpha_{n}}
               (\, _{\tau}D_{t}^{\sigma_{n-1}}G(t,\tau)) f(\tau) d\tau
                           \nonumber \\
        & & \hspace*{5em}
           + \lim_{\tau \rightarrow t-0}
               \, _{\tau}D_{t}^{\alpha_{n} -1}
                (\, _{\tau}D_{t}^{\sigma_{n-1}} G(t,\tau)) f(\tau)
                           \nonumber \\
        & = & \int_{0}^{t}\, _{\tau}D_{t}^{\sigma_{n}}G(t,\tau) f(\tau) d\tau
              + f(t)
\end{eqnarray}

Multiplying these equations by the corresponding coefficients
and summarizing, we obtain
\begin{equation}
      _{0}{\cal L}_{t}y(t) = \int_{0}^{t} \, _{\tau}{\cal L}_{t}
                              G(t, \tau) f(\tau) d\tau + f(t) = f(t),
\end{equation}
which completes this proof.

\vspace{1.5ex}

\underline{{\bf 2.}}
For fractional differential equations with constant coefficients
$G(t, \tau) \equiv G(t - \tau)$.
This is obvious because in such a case
Green's function can be obtained by Laplace transform method.

\vspace{1.5ex}

The type (standard or sequential) of the equation
is not important for determining Green's function,
because due to condition (b) in the Green's function definition
all non-integral addends vanish.

\vspace{1.5ex}

\underline{{\bf 3.}}
Appropriate derivatives
of Green's function $G(x,\tau)$ form a set of linearly independent
solutions of a homogeneous $(f(t) \equiv 0)$ equation (\ref{IVP})
(for a simple illustration, see Examples \ref{Samko} and \ref{Example-3a}).

\vspace{1.5ex}

Let us demonstrate this for the case of the linear fractional differential
equations with constant coefficients, which are the main
subject for study in this work and for which we have
$G(t, \tau) \equiv G(t-\tau)$.

Let us take $0< \lambda < \sigma_n$. First, the function
\begin{equation}
      y_{\lambda}(t) = \,\, _{0}D_{t}^{\lambda} G(t)
\end{equation}
is a solution of the corresponding homogeneous equation. Indeed,
\begin{equation}
      _{0}{\cal L}_{t}y_{\lambda}(t) =
      \,\, _{0}{\cal L}_{t} (\,\, _{0}D_{t}^{\lambda}G(t)) =
      \,\, _{0}D_{t}^{\lambda}(\,\,_{0}{\cal L}_{t} G(t)) = 0.
\end{equation}
We used here that
$ _{0}{\cal L}_{t} \, _{0}D_{t}^{\lambda} = \,\,
        _{0}D_{t}^{\lambda} \, _{0}{\cal L}_{t}$,
which follows from condition (c) in the definition
of the fractional Green's function.

Second,
\begin{equation}
      \left.
      _{0}D_{t}^{\sigma_n - \lambda -1} y_{\lambda} (t)
      \right|_{t=0}
      = 1.
\end{equation}

In fact,
\begin{equation}
      \left.
      _{0}D_{t}^{\sigma_n - \lambda -1} y_{\lambda} (t)
      \right|_{t=0}
      =   \left. \,
          _{0}D_{t}^{\sigma_n - \lambda -1}
          \left(
                \,\,_{0}D_{t}^{\lambda} G(t)
          \right)
          \right|_{t=0}
       =  \left.
           \,\,_{0}D_{t}^{\sigma_n -1} G(t)
           \right|_{t=0}
       = 1.
\end{equation}

We first used here the relationship
\begin{equation}
      _{0}D_{t}^{\sigma_n -\lambda -1} \,\,
      _{0}D_{t}^{\lambda} G(t)
        \equiv
      \,\, _{0}D_{t}^{(\sigma_n -\lambda -1)+\lambda} G(t),
\end{equation}
which follows from condition (c) of the definition
of Green's function, and then condition (b).

We see that having the fractional Green's function of
equation (\ref{IVP-G}), we can determine particular solutions
of the corresponding homogeneous equation,
which are necessary for satisfying inhomogeneous initial conditions.

Therefore, solution of linear fractional differential equations
with constant coefficients reduces to finding
the fractional Green's function. After that, we can
immediately write the solution
of the inhomogeneous equation satisfying given
inhomogeneous initial conditions.

This solution has the form
\begin{equation}
      y(t) = \sum_{k=1}^{n}b_k \psi_k (t)
              + \int_{0}^{t} G(t - \tau) f(\tau) d\tau,
\end{equation}
\begin{equation}
      b_k= \left.
                \, _{0}{\cal D}_{t}^{\sigma_k -1}y(t)
           \right|_{t=0}
\end{equation}
\begin{equation}
      \psi_k (t) = \,\, _{0}{\cal D}_{t}^{\sigma_n -\sigma_k}G(t),
      \qquad
      _{0}{\cal D}_{t}^{\sigma_n -\sigma_k}
       \equiv
      \,  _{a}D_{t}^{\alpha_{n}}
      \, _{a}D_{t}^{\alpha_{n-1}} \ldots
      \,_{a}D_{t}^{\alpha_{k+1}}
\end{equation}

Because of this, in the following sections we find some
explicit expressions for fractional Green's functions,
including a general linear fractional differential equation.

\section{Fractional Green's function  for the one-term \protect \\
         fractional  differential equation}

The fractional Green's function $G_{1}(t)$
for the one-term fractional-order differential equation
with constant coefficients
\begin{equation} \label{One-term-equation}
     a\,\, _{0}D_{t}^{\alpha}y(t) = f(t),
\end{equation}
where the derivative can be either
"classic" (i.e., considered in the book by Oldham and Spanier)
or "sequential" (Miller and Ross),
is found by the inverse Laplace transform
of the following expression:
\begin{equation}
      g_1 (p) = \frac{1}{ap^\alpha}.
\end{equation}

The inverse Laplace transform then gives
\begin{equation}
      G_1 (t) = \frac{1}{a} \,\, \frac{t^{\alpha -1}}{\Gamma (\alpha)}.
\end{equation}

The solution of equation (\ref{One-term-equation}) under homogeneous
initial conditions is
\begin{equation} \label{G1-solution}
      y(t) = \frac{1}{a \Gamma (\alpha)}
             \int_{0}^{t}
             \frac{f(\tau) d\tau}{(t-\tau)^{1-\alpha}}
           = \frac{1}{a} \,\, _{0}D_{t}^{-\alpha}f(t).
\end{equation}

Using \cite[lemma 3.3]{Pitcher-Sewell} , we can easily verify that
expression (\ref{G1-solution}) gives the solution of equation
(\ref{One-term-equation}), if $f(x)$ is continuous in $[0, \infty)$.

\section{Fractional Green's function
         for the two-term \protect \\  fractional
         differential equation}

The fractional Green's function $G_{2}(t)$
for the two-term fractional-order differential equation
with constant coefficients
\begin{equation} \label{Two-term-equation}
     a\,\, _{0}D_{t}^{\alpha}y(t) + by(t) = f(t),
\end{equation}
where the derivative can be either
"classic" (i.e., considered in the book by Oldham and Spanier)
or "sequential" (Miller and Ross),
is found by the inverse Laplace transform
of the following expression:
\begin{equation}
      g_2 (p) = \frac{1}{ap^\alpha +b}
              = \frac{1}{a} \,\, \frac{1}{p^\alpha + \frac{b}{a}},
\end{equation}
which leads to
\begin{equation}\label{G2}
      G_2 (t) = \frac{1}{a} t^{\alpha -1}
                    E_{\alpha, \alpha}(-\frac{b}{a}t^\alpha).
\end{equation}
For instance, function $G_2 (t)$ plays the key role in solution
of Example \ref{Oldham-Spanier-1} and \ref{Samko}.

Taking in (\ref{G2}) $b=0$ and using the definition
of the Mittag-Leffler function
(\ref{ML-Definition}), we obtain Green's function
$G_1 (t)$ for the one-term equation.

\section{Fractional Green's function
         for the three-term \protect \\  fractional
          differential equation}

The fractional Green's function $G_{3}(t)$
for the three-term fractional-order differential equation
with constant coefficients
\begin{equation} \label{Three-term-equation}
     a\,\, _{0}D_{t}^{\beta}y(t) +
     b\,\, _{0}D_{t}^{\alpha}y(t) +
     c\,\, y(t) = f(t),
\end{equation}
where the derivatives can be either
"classic" (i.e., considered in the book by Oldham and Spanier)
or "sequential" (Miller and Ross),
is found by the inverse Laplace transform
of the following expression:

\begin{equation}
      g_{3}(p) = \frac{1}{ap^{\beta} + bp^{\alpha} + c}
\end{equation}

Assuming $\beta > \alpha$, we can write $g_{3}(p)$ in the form
\begin{equation}
  g_{3}(p)  =  \frac{1}{c} \,\,
             \frac{cp^{-\alpha}}{ap^{\beta - \alpha} + b} \,\,
             \frac{1}{1+ \frac{\mbox{$cp^{-\alpha}$}}
                              {\mbox{$ap^{\beta - \alpha} +b$}}}
        =  \frac{1}{c}
             \sum_{k=0}^{\infty}
             (-1)^{k}
             \left( \frac{c}{a} \right)^{k+1}
             \frac{p^{-\alpha k -\alpha}}
                     {\left( p^{\beta - \alpha} +
                         \frac{\mbox{$b$}}{\mbox{$a$}} \right)^{k+1} }
\end{equation}
The term-by-term inversion, based
on the general expansion theorem for the Laplace transform
given in \cite[\S 22]{Doetsch},
using (\ref{Laplace-Transform}) produces
\begin{equation} \label{My-Solution}
      G_{3}(t) = \frac{1}{a}
             \sum_{k=0}^{\infty}
             \frac{(-1)^k}{k!}
             \left( \frac{c}{a}  \right)^{k}
             t^{\beta (k+1) -1}
             E_{\beta - \alpha, \beta + \alpha k}^{(k)}
               (-\frac{b}{a} t^{\beta - \alpha}),
\end{equation}
where $E_{\lambda, \mu}(z)$ is the
Mittag-Leffler function in two parameters,
\begin{equation}\label{E-Derivative}
E_{\lambda,\mu}^{(k)}(y) \equiv \frac{d^{k}}{dy^{k}}E_{\lambda ,\mu}(y) =
\sum_{j=0}^{\infty} \frac{(j+k)! \,\, y^{j}}
                         {j! \,\, \Gamma (\lambda j + \lambda k + \mu)},
\hspace{3em}
(k = 0, 1, 2, ...)
\end{equation}

We assume in this solution that $a \neq 0$, because otherwise we have
the two-term equation (\ref{Two-term-equation}).
We can also assume $c \neq 0$, because for $c=0$
$$
   g_3 (p) = \frac{1}{ap^\beta +bp^\alpha}
           = \frac{p^{-\alpha}}{ap^{\beta- \alpha} + b},
$$
and the Laplace inversion can be done in the same way
as in the case of the two-term equation.

Two special cases of equation (\ref{Three-term-equation})
were considered by Bagley and Torvik \cite{Bagley-Torvik}
(for $\beta = 2$ and $\alpha = 3/2$)
and by Caputo \cite{Caputo-ED} (for $\beta =2$ and $0 < \alpha < 1$).
It is easy to show that our solution (\ref{My-Solution})
contains Caputo's solution as a particular case.

Indeed, substituting (\ref{E-Derivative}) into (\ref{My-Solution}) and
changing the order of summation, we obtain:
\begin{eqnarray}
\hspace{-1em}
  G_{3}(t) & = & \frac{1}{c}
             \sum_{k=0}^{\infty}
             (-1)^{k}
             \left( \frac{c}{a}  \right)^{k+1}
             \sum_{j=0}^{\infty}
             (-1)^{j}
             \left( \frac{b}{a} \right)^{j}
             \frac{(j+k)! \,\, t^{\beta (j+k) +\beta -1 -\alpha j}}
                  {k! \,\, j! \,\, \Gamma (\beta (j+k+1) - \alpha j)}
                            \nonumber \\
  \hspace{-1em}
       & = & \frac{1}{c}
             \sum_{j=0}^{\infty}
             \left( \frac{-b}{a}  \right)^{j}
             \sum_{k=0}^{\infty}
             (-1)^{k}
             \left( \frac{c}{a} \right)^{k+1}
             \frac{(j+k)! \,\, t^{\beta (j+k) +\beta -1 -\alpha j}}
                  {k! \,\, j! \,\, \Gamma (\beta (j+k+1) - \alpha j)}
\end{eqnarray}
For $\beta =2$ this expression is identical with the expression
obtained by Caputo \cite[formula (2.27)]{Caputo-ED}.

\section{Fractional Green's function
         for the four-term \protect \\ fractional
         differential equation}

The fractional Green's function $G_{4}(t)$
for the four-term fractional-order differential equation
with constant coefficients
\begin{equation} \label{Four-term-equation}
     a\,\, _{0}D_{t}^{\gamma}y(t) +
     b\,\, _{0}D_{t}^{\beta}y(t) +
     c\,\, _{0}D_{t}^{\alpha}y(t) +
     d\,\, y(t) = f(t),
\end{equation}
where the derivatives can be, as in the previous section, either
"classic"
or "sequential",
is found by the inverse Laplace transform
of the following expression:
\begin{equation}
      g_{4}(p) = \frac{1}{ap^{\gamma} + bp^{\beta} + cp^{\alpha} + d}
\end{equation}

Assuming $\gamma > \beta > \alpha$, we can write $g(p)$ in the form
\begin{eqnarray}
  g_{4}(p) & = & \frac{1}{ap^{\gamma}+bp^{\beta}} \,\,
             \frac{1}{1+ \frac{\mbox{$cp^{\alpha}$}+\mbox{$d$}}
                              {\mbox{$ap^{\gamma}+bp^{\beta}$}} }
                                             \nonumber \\
       & = & \frac{a^{-1}p^{-\beta}}
                     {p^{\gamma -\beta}+ a^{-1}b} \,\,
             \frac{1}{1+ \frac{\mbox{$a^{-1}cp^{\alpha -\beta}$}
                                     +\mbox{$a^{-1}dp^{-\beta}$}}
                              {\mbox{$p^{\gamma -\beta}+a^{-1}b$}} }
                                             \nonumber \\
      & = & \sum_{m=0}^{\infty}
            (-1)^{m}
            \frac{a^{-1}p^{-\beta}}{(p^{\gamma -\beta} + a^{-1}b)^{m+1}}
            \left(
               \frac{c}{a}p^{\alpha -\beta} + \frac{d}{a}p^{-\beta}
            \right)^{m}
                                             \nonumber \\
      & = & \sum_{m=0}^{\infty}
            (-1)^{m}
            \frac{a^{-1}p^{-\beta}}{(p^{\gamma -\beta} + a^{-1}b)^{m+1}}
            \sum_{k=0}^{m}
            {m \choose k}
             \frac{c^{k}d^{m-k}}{a^m}
             p^{\alpha k - \beta m}
                                        \nonumber \\
      & = & \frac{1}{a}
            \sum_{m=0}^{\infty}
              (-1)^{m} \left( \frac{d}{a} \right)^m
            \sum_{k=0}^{m}
              {m \choose k}
               \left( \frac{c}{d} \right)^k
              \frac{p^{\alpha k - \beta m - \beta}}
                   {(p^{\gamma - \beta} + a^{-1}b)^{m+1}}
\end{eqnarray}

The term-by term inversion, based
on the general expansion theorem for the Laplace transform
given in \cite[\S 22]{Doetsch},
using (\ref{Laplace-Transform}) gives the final
expression for the fractional Green's function
for equation (\ref{Four-term-equation}):
\begin{equation}
  G_{4}(t) \! = \! \frac{1}{a}
            \sum_{m=0}^{\infty}\frac{1}{m!}
              \left( \frac{-d}{a} \right)^m
            \sum_{k=0}^{m}
              {m \choose k}
               \left( \frac{c}{d} \right)^k
             t^{\gamma (m+1) - \alpha k - 1}
             E_{\gamma -\beta, \gamma +\beta m -\alpha k}^{(m)}
                   (-\frac{b}{a}t^{\gamma -\beta})
\end{equation}

We assumed in this solution $a \neq 0$,
because in the opposite case we have
the three-term equation (\ref{Three-term-equation}).
We can also assume $d \neq 0$, because in the case of $d = 0$,
after writing
\begin{equation}
  g_{4}(p) = \frac{p^{-\alpha}}{ap^{\gamma -\alpha} + bp^{\beta -\alpha} + c}
\end{equation}
the Laplace inversion can be done in the same way as
in the case of the three-term equation.

\section{Fractional Green's function
         for the general \protect \\  linear fractional
         differential equation }

The above results can be essentially generalized.

The fractional Green's function $G_{n}(t)$
for the $n$\/-term fractional-order differential equation
with constant coefficients
\begin{equation} \label{n-term-equation}
     a_{n}\, D^{\beta_{n}}y(t) +
     a_{n-1}\, D^{\beta_{n-1}}y(t) + \ldots +
     a_{1}\, D^{\beta_{1}}y(t) +
     a_{0}\, D^{\beta_{0}}y(t) = f(t),
\end{equation}
where derivatives $D^{\alpha} \equiv \,\, _{0}D_{t}^{\alpha}$ can be, as in the previous sections, either
"classic" or "sequential",
is found by the inverse Laplace transform
of the following expression:
\begin{equation}
      g_{n}(p) = \frac{1}{a_{n}p^{\beta_{n}} + a_{n-1}p^{\beta_{n-1}} +\ldots
                  + a_{1}p^{\beta_{1}} + a_{0}p^{\beta_{0}}}
\end{equation}

Let us assume
$\beta_{n} > \beta_{n-1} > \ldots > \beta_{1} > \beta_{0}$ and
write $g_{n}(p)$ in the form:
\begin{eqnarray}
  g_{n}(p) & = & \frac{1}{a_{n}p^{\beta_{n}}+a_{n-1}p^{\beta_{n-1}}}
               \,
               \frac{1}
                    {1+ \frac{\mbox{$\sum\limits_{k=0}^{n-2}a_{k}p^{\beta_{k}}$}}
                         {\mbox{$a_{n}p^{\beta_{n}}+a_{n-1}p^{\beta_{n-1}}$}}}
                                \nonumber \\
   & = & \frac{a_{n}^{-1}p^{-\beta_{n-1}}}{p^{\beta_{n}-\beta_{n-1}}
                     +\frac{\mbox{$a_{n-1}$}}{\mbox{$a_n$}}}
               \,
               \frac{1}
                    {1 + \frac{\mbox{$a_{n}^{-1}p^{-\beta_{n-1}}$}
                              \mbox{$\sum\limits_{k=0}^{n-2}a_{k}p^{\beta_{k}}$}}
                            {\mbox{$p^{\beta_{n}-\beta_{n-1}}
                                 + \frac{\mbox{$a_{n-1}$}}{\mbox{$a_n$}}
                             $}}}
                                      \nonumber \\
   & = & \sum_{m=0}^{\infty}
         \frac{\mbox{$(-1)^{m} a_{n}^{-1} p^{-\beta_{n-1}}$}}
              {\left(
                      \mbox{$p^{\beta_{n} -\beta_{n-1}}
                      + \frac{\mbox{$a_{n-1}$}}{\mbox{$a_n$}}$}
               \right)^{m+1}}
          \left(
             \sum_{k=0}^{n-2}
             \left(
                \frac{\mbox{$a_k$}}{\mbox{$a_n$}}
             \right)
             p^{\beta_{k}-\beta_{n-1}}
          \right)^m
                        \nonumber \\
     & = & \sum_{m=0}^{\infty}
          \frac{(-1)^{m} a_{n}^{-1} p^{-\beta_{n-1}} }
               {\left(
                    p^{\beta_n -\beta_{n-1}} + \frac{\mbox{$a_{n-1}$}}
                                                   {\mbox{$a_n$}}
                \right)^{m+1} }
                \hspace*{10em} \mbox{(continued)}
                                   \nonumber \\
& &  \hspace*{1em}   \sum_{{ k_0+k_1+\ldots +k_{n-2} =m
                  \atop
             k_0 \geq 0; \ldots\, , k_{n-2} \geq 0  }}
          (m; k_0, k_1, \ldots \,, k_{n-2})
          \prod_{i=0}^{n-2}
     \left(
          \frac{a_i}{a_n}
     \right)^{k_i}
      p^{(\beta_i-\beta_{n-1})k_i}
                                   \nonumber \\
& = & \frac{1}{a_n} \sum_{m=0}^{\infty} (-1)^m
     \sum_{{ k_0+k_1+\ldots +k_{n-2} =m
                  \atop
             k_0 \geq 0; \ldots\, , k_{n-2} \geq 0
            }}
     (m; k_0, k_1, \ldots \,, k_{n-2}) \hspace*{1em} (\mbox{continued})
                                   \nonumber \\
& &  \hspace*{8em}
     \prod_{i=0}^{n-2}
     \left(
          \frac{a_i}{a_n}
     \right)^{k_i}
     \frac{p^{-\beta_{n-1} +\sum\limits_{i=0}^{n-2}(\beta_i-\beta_{n-1})k_i }}
          {\left(
                    p^{\beta_n -\beta_{n-1}} + \frac{\mbox{$a_{n-1}$}}
                                                   {\mbox{$a_n$}}
                \right)^{m+1}}
\end{eqnarray}
where $((m; k_0, k_1, \ldots \,, k_{n-2}))$ are the multinomial coefficients
\cite[chapter 24]{Abramowitz-Stegun}

The term-by term inversion, based
on the general expansion theorem for the Laplace transform
given in \cite[\S 22]{Doetsch},
using (\ref{Laplace-Transform}) gives the final
expression for the fractional Green's function
for equation (\ref{n-term-equation}):
\begin{eqnarray}
     G_{n}(t) & = & \frac{1}{a_n} \sum_{m=0}^{\infty} \frac{(-1)^m}{m!}
     \sum_{{ k_0+k_1+\ldots +k_{n-2} =m
                  \atop
             k_0 \geq 0; \ldots\, , k_{n-2} \geq 0
            }}
     (m; k_0, k_1, \ldots \,, k_{n-2}) \hspace*{1em} (\mbox{continued})
                                   \nonumber \\
& &  \prod_{i=0}^{n-2}
     \left(
          \frac{a_i}{a_n}
     \right)^{k_i}
     t^{(\beta_n -\beta_{n-1})m +\beta_n
                      +\sum_{j=0}^{n-2}(\beta_{n-1}-\beta_j)k_j -1}
                                        \hspace*{2em} (\mbox{continued})
                                   \nonumber \\
& &  E^{(m)}_{\beta_n-\beta_{n-1}, +\beta_n
                      +\sum_{j=0}^{n-2}(\beta_{n-1}-\beta_j)k_j }
         \left(
               -\frac{a_{n-1}}{a_n}t^{\beta_n-\beta_{n-1}}
         \right)
\end{eqnarray}

\end{sloppypar}

\newpage
\thispagestyle{empty}
\vspace*{17cm}
\noindent
N\'azov: The Laplace Transform Method for Linear Differential \\
\hspace*{3.5em}       Equations of the Fractional Order \\
Autor:  RNDr. Igor Podlubny, CSc. \\
Zodp. redaktor: RNDr. P.Samuely, CSc.\\
Vydavate\softl{}: \'Ustav experiment\'alnej fyziky SAV, Ko\v{s}ice \\
Redakcia: \'UEF SAV, Watsonova 3, 04001 Ko\v{s}ice, Slovensk\'a republika \\
Po\v{c}et str\'an: 32 \\
N\'aklad: 120  \\
Rok vydania: 1994 \\
Tla\v{c}: OLYMPIA s.r.o.,  M\'anesova 23, 040 01 Ko\v{s}ice


\begin{thebibliography}{99} \addcontentsline{toc}{chapter}{Bibliography}
\bibitem{OS} K.B.Oldham and J.Spanier.
                {\em The fractional calculus}.
                Academic Press, New York, 1974.
\bibitem{Miller-Ross} K.S.Miller and B.Ross.
                {\em An introduction
                to the fractional calculus and fractional
                differential equations}.
                John Wiley \& Sons. Inc., New York, 1993.
\bibitem{Bagley-Calico} R.L.Bagley and R.A.Calico.
                {\em J.Guidance}, vol.14, No.2, 1991,
                pp.304-311.
\bibitem{Caputo-Mainardi} M.Caputo and F.Mainardi.
                {\em Rivista del Nuevo Cimento (Serie II)}, vol.1, No.2,
                1971, pp.161-198.
\bibitem{SKM} S.G.Samko, A.A.Kilbas and O.I.Maritchev.
                {\em Integrals and derivatives of the fractional order
                and some of their applications}.
                Nauka i Tekhnika,
                Minsk, 1987 (in Russian).
\bibitem{Schneider-Wyss} W.R.Schneider and W.Wyss.
               {\em J.Math.Phys.}, {\bf 30}, 1989,
               pp. 134-144.
\bibitem{GN} W.G.Gl\"{o}ckle and T.F.Nonnenmacher.
                {\em Macromolecules}, vol.24, 1991, pp. 6426-6436.
\bibitem{Fox} C.Fox.
               {\em Trans.Am.Math.Soc.}, vol.98, 1961,
               pp.395-429.
\bibitem{HTF-1} A.Erd\'elyi (ed.). {\em Higher transcendental functions}, vol.1.
                McGraw-Hill, New York, 1955.
\bibitem{HTF-3} A.Erd\'elyi (ed.). {\em Higher transcendental functions}, vol.3.
                McGraw-Hill, New York, 1955.

\bibitem{Rabotnov} Yu.N.Rabotnov. {\em Elements} {\em of} {\em hereditary}
                                  {\em solids} {\em mechanics}.
               Na\-uka, Mos\-cow, 1977 (in Russian).
\bibitem{Oustaloup}  A.Oustaloup.
               {\em Proceedings of 12th IMACS World Congress},
               Paris, July 18-22, 1988, vol.3,  pp. 203-208.
\bibitem{Nigmatullin-Transfer} R.R.Nigmatullin.
               {\em Physica Status Solidi (B)},
               {\bf 133}, 1986, pp.425-430.

\bibitem{Westerlund} S.Westerlund.
               {\em Physica scripta}, vol.43, 1991, pp.174-179.
\bibitem{Mainardi-WASCOM} F.Mainardi. {\em On the initial value problem
               for the fractional diffusion-wave equation}. In the book:
               S.Rionero and T.A. Ruggeri (eds.) Proceedings of VII WASCOM,
               World Scientific, Singapore, 1994 (pre-print)
\bibitem{Friedrich} Ch.Friedrich.
               {\em J.Non-Newtonian Fluid Mech.}, vol. 46, 1993, pp.307-314.
\bibitem{Friedrich-1991} Ch.Friedrich.
                {\em Rheologica Acta}, 30, 1991, pp.151-158.

\bibitem{Caputo-ED} M.Caputo. {\em Elasticit\`a e dissipazione}.
               Zanichelli, Bologna, 1969.

\bibitem{TIT-I} A.Erd\'elyi (ed.). {\em Tables of integral transforms},
                vol. 1, McGraw-Hill, 1954.
\bibitem{Doetsch}G.Doetsch. {\em Anleitung} {\em zum praktischen}
                 {\em gebrauch}
                 {\em der} {\em Laplace-trans\-form\-ation}.
                 Ol\-den\-bourg, M\" unchen, 1956.
                 (Russian translation: Fizmatgiz, Moscow, 1958).
\bibitem{Bagley-Torvik} R.L.Bagley and P.J.Torvik.
               {\em J. Appl. Mech.}, vol.51, 1984, pp.294-298.
\bibitem{Abramowitz-Stegun} M.Abramowitz and I.A.Stegun.
                {\em Handbook of Mathematical Functions.}
                Nat. Bureau of Standards, Appl. Math. Series, {\bf 55}, 1964
                (Russian translation: Nauka, Moscow, 1979.)
\bibitem{Agarwal}R.P.Agarwal. A propos d'une note de M.Pierre Humbert.
                 {\em C.R. Seances Acad. Sci.},
                 1953, vol. 236, No.21, pp.2031-2032.
\bibitem{Humbert-Agarwal}P.Humbert et R.P.Agarwal.
                 Sur la fonction de Mittag-Leffler et
                 quel\-ques\--unes de ses g\'en\'eralisations.
                 {\em Bulletin des Sciences Mathematiques},
                 1953, vol. 77, No.10, pp.180-185.
\bibitem{Plotnikov}Yu.I.Plotnikov. {\em Steady-state vibrations of plane
                and axesym\-met\-ric \\ stamps on a viscoelastic foundation}.
                PhD thesis, Moscow, 1979 (in Russian).
\bibitem{Tseytlin} A.I.Tseytlin. {\em Applied methods of solution of boundary
                value problems in civil engineering}.
                Stroyizdat, Moscow, 1984 (in Russian).
\bibitem{Naimark}M.A.Naimark. {\em Linear differential operators}.
                Nauka, Moscow, 1969 (in Russian).
\bibitem{Pitcher-Sewell} E.Pitcher and W.E.Sewell.
                {\em Bull. Amer. Math. Soc.},
                vol.44, No. 2, 1938, pp.100-107.
\bibitem{Wright-33} E.M.Wright. {\em J. London Math.Soc.}, 1933,
                vol.8, pp.71-79.
                {\begin{center} *** *** *** \end{center}}
\bibitem{Caputo-1967} M.Caputo. {\em Geophys.J.R.Astr.Soc.}, 1967, vol. 13,
                      pp.529-539.
\bibitem{Caputo-1993} M.Caputo. {\em Rend. Fis. Acc.Lincei}, s. 9, v.4, pp.89-98.
\bibitem{Caputo-1989} M.Caputo. {J.Acoust.Soc.Am.}, vol. 86, No. 5 (Nov 1989),
                      pp.1984-1987.
\bibitem{Caputo-1990} F.Bella, P.F.Biagi, M.Caputo, G. Della Monica,
                      A.Ermini, P.Manjgaladze, V.Sgrigna and D.Zilpimiani.
                      {\em Tectonophysics}, vol. 179, 1990, pp. 131-139.
\bibitem{Caputo-Lectures} M.Caputo. {\em Lectures on seismology and
                      rheological tectonics}. Univ. deglistudi di Roma
                      ''La Sapienza'', 1992-1993.
%
\bibitem{Macris-1993} N.Macris, G.F.Dargush and M.C.Constantinou.
                      {\em ASCE J. Eng. Mech.}, 1993, vol. 119, No. 8,
                      pp. 1663-1679.
\bibitem{Macris-1991} N.Macris and M.C.Constantinou. {\em J.Fluid Mech.},
                      1991, vol. 225, pp. 631-653.
%
\bibitem{Ochmann-1993} M.Ochmann and S.Makarov. {\em J.Amer.Acoust.Soc.},
                      1993, vol. 94, No. 6 (Dec 1993), pp. 3392-3399.
%
\bibitem{Friedrich-1992} Ch.Friedrich. {\em Philosophical Magazine Letters},
                        1992, vol. 66, No. 6, pp.287-292.

\bibitem{Giona-1992a} M.Giona and H.E.Roman.{\em The Chemical Engineering Journal},
                        1992, vol. 49, pp.1-10.
\bibitem{Giona-1992b} M.Giona, S.Gerbelli and H.E.Roman. {Physica A},
                        1992, vol. 191, pp.449-453.

\bibitem{Bagley-1989} R.L.Bagley. {\em AIAA Journal}, 1989, vol. 27, No. 10
                         (Oct 1989), pp.1412-1417.

\bibitem{Koeller-1984} R.C.Koeller. {\em Trans. ASME - J.Appl.Mech.},
                        1984, vol. 51, pp. 299-307.

\bibitem{Nonnenmacher-1989} T.F.Nonnenmacher and D.J.F.Nonnenmacher.
                       {\em Acta Physica Hungarica}, 1989, vol. 66,
                       pp.145-154.
\bibitem{Nonnenmacher-1991} T.F.Nonnenmacher and W.G.Gl\"ockle.
                        {\em Phylosophical Magazine Letters}, 1991, vol. 64,
                        No. 2, pp. 89-93.
\bibitem{Nonnenmacher-1990} T.F.Nonnenmacher.
                       {\em J. of Physics A: Math. and Gen.},
                       1990, vol. 23, pp. L697-L700.
\bibitem{Nonnenmacher-1991-L} T.F.Nonnenmacher, in:
                       {\em Lect. Notes in Physics, vol. 381},
                       Springer-Verlag, Berlin, 1991,
                       pp. 309-320.

\bibitem{He-1987} X.-F.He {\em Solid State Comm.}, 1987, vol. 61, No. 1,
                        pp. 53-55.
\bibitem{He-1993} D.Mo, Y.Y.Lin, J.H.Tan, Z.X.Yu, G.Z.Zhou,
                  K.C.Gong, G.P.Zhang, X.-F.He
                  {\em Thin Solid Films}, 1993, vol. 234, pp. 468-470.

\bibitem{Le Mehaute-1983} A.Le~Mehaute and G.Crepy.
                        {\em Solid State Ionics}, 1983, No. 9 and 10,
                        pp. 17-30.

\bibitem{Nigmatullin-1984-a} R.R.Nigmatullin. {\em Phys. Sta.Sol (b)},
                        1984, vol. 123, pp. 739-745.
\bibitem{Nigmatullin-1984-b} R.R.Nigmatullin. {\em Phys. Sta.Sol (b)},
                        1984, vol. 124, pp. 389-393.
\bibitem{Nigmatullin-1992} R.R.Nigmatullin. {\em Soviet J. Theor.and Math. Phys.},
                        1992, vol. 90, No.3, pp. 354-367.

\bibitem{Kaloyanov-1992} G.Kaloyanov and J.M.Dimitrova.
                        {\em Izv. Vysshykh Utchebnykh Zav. - Elektromekhanika},
                        1992, No.2, pp. 65-72  (in Russian)

\bibitem{Axtell-1990} M.Axtell and M.E.Bise, in:
                      {\em Proc. of the IEEE 1990 Nat. Aerospace and Electronics
                      Conf., New York}, 1990,
                      pp. 563-566.


\end{thebibliography}
\end{document}